\documentclass{article}
\usepackage{amsmath, amsthm}

\usepackage{amssymb,amscd,amsfonts,euscript}

\newtheoremstyle{theorem}
  {10pt}		  
  {10pt}  
  {\sl}  
  {\parindent}     
  {\bf}  
  {. }    
  { }    
  {}     
\theoremstyle{theorem}
\newtheorem{theorem}{Theorem}

\newtheorem{proposition}[theorem]{Proposition}
\newtheorem{lemma}[theorem]{Lemma}

\newtheoremstyle{defi}
  {10pt}		  
  {10pt}  
  {\rm}  
  {\parindent}     
  {\bf}  
  {. }    
  { }    
  {}     
\theoremstyle{defi}
\newtheorem{definition}[theorem]{Definition}
\newtheorem{remark}[theorem]{Remark}

\def\PM{{\mathbb{P}}}
\def\EM{{\mathbb{E}}}
\def\NM{{\mathbb{N}}}
\def \brond{\mathcal{B}}
\def \frond{\mathcal{F}}
\def \grond{\mathcal{G}}
\def \trond{\mathcal{T}}
\def \card {\mathop{\hbox{\rm Card}}\nolimits}

\def \guillegauche{\/\textquotedblleft\/}
\def \guilledroit{\/\textquotedblright\/}
\def \va#1{{\bf #1}}
\def \implique {\mathop{\Longrightarrow}\limits}
\renewcommand\leq{\leqslant}
\renewcommand\geq{\geqslant}
\newcounter{c}
\def \cpp {\addtocounter{c}{1}}
\def \con {c_{\thec}}

\begin{document}

\title{A quasimartingale characterization of $p$-stable type Banach spaces.}

\author{Florian HECHNER $^1$ and Bernard HEINKEL$^2$\\
$^1$Corresponding author\\IRMA, University of Strasbourg and CNRS,\\ 7 rue Ren\'e-Descartes, 67084 STRASBOURG CEDEX\\
hechner@math.unistra.fr\\[2pt]
$^2$IRMA, University of Strasbourg and CNRS,\\ 7 rue Ren\'e-Descartes, 67084 STRASBOURG CEDEX\\
heinkel@math.unistra.fr}

\maketitle

\begin{abstract}
We characterize Banach spaces $\brond$ of stable-type $p$ ($1<p<2$) by the property that for every sequence $({\bf X}_i)$ of ${\brond}$-valued random variables, independent, centered and fulfilling some integrability assumption, the sequence $\left(\frac{{\bf X}_1+\cdots+{\bf X}_n}{n^{1/p}}\right)$ is a quasimartingale.


\end{abstract}

\section{Introduction}

~\guillegauche How to characterize the regularity of a Banach space $(\brond,\|\cdot\|)$ by the fact that a kind of classical strong law of large numbers (SLLN) holds for $\brond$-valued random variables (r.v.)?\guilledroit\ is a well know problem. Two kinds of regular Banach spaces -- the spaces of Rademacher type $p$ and the spaces of stable type $p$ -- have been nicely characterized in that way (see chapter 9 in \cite{Ledoux1}).

Here our goal will be to show that the SLLN characterization of spaces of stable type $p$ ($1<p<2$) can be made more precise in supposing that the \guillegauche normalized sums\guilledroit\ obeying the SLLN have even a quasimartingale behaviour.

To begin with, we recall some definitions.

In the sequel, $(\brond,\|\cdot\|)$ will be a real separable Banach space, equipped with its Borel $\sigma$-field $B$. A $\brond$-valued r.v. $\va X$ is a measurable function defined on a probability space $(\Omega,\trond,\PM)$ with values in $(\brond,B)$. Such a r.v. is said to be (strongly) integrable if $\EM\|\va X\|<+\infty$ and $\forall f\in\brond^\prime,\ \EM f(\va X)=0$; this is denoted by $\EM(\va X)=0$.

Let $p\geq 1$ be given. The weak-$\ell_p$ norm of a sequence $a:=(a_1,\ldots,a_n)$ of real numbers is defined as follows~:
$$
\|a\|_{p,\infty}:=\sup\limits_{t>0}{\left(t^p\card(i:|a_i|>t)\right)}^\frac{1}{p}=\sup\limits_{k=1}^na_k^*k^{1/p},
$$
where $(a_1^*,\ldots,a_n^*)$ denotes the non-increasing rearrangement of the sequence $(|a_1|,\ldots,|a_n|)$.

Let now $(\varepsilon_k)$ be a sequence of independent Rademacher random variables $\left(\text{that is}\ \PM(\varepsilon_k=1)=\PM(\varepsilon_k=-1)=\frac{1}{2}\right)$.

Rademacher type $p$ spaces and stable type $p$ spaces are defined as follows~:

\begin{definition}\label{defi1}~
\begin{enumerate}
\item Let $1<p<2$. The Banach space $(\brond,\|\cdot\|)$ is of Rademacher type $p$ if there exists a constant $c(p)>0$ such that for every finite sequence $(x_i)$ in $\brond$~:
$$
{\left(\EM\left\|\sum\limits_{i=1}^n\varepsilon_i x_i\right\|^p\right)}^{\frac{1}{p}}\leq c(p){\left(\sum\limits_{i=1}^n\|x_i\|^p\right)}^{\frac{1}{p}}.
$$
\item Let $1\leq p\leq 2$. The Banach space $(\brond,\|\cdot\|)$ is of stable type $p$ if there exists a constant $C(p)>0$ such that for every finite sequence $(x_i)$ in $\brond$~:
$$
\EM\left\|\sum\limits_{i=1}^n\varepsilon_i x_i\right\|\leq C(p)\left\|\left(\|x_1\|,\ldots,\|x_n\|\right)\right\|_{p,\infty}.
$$
\end{enumerate}
\end{definition}

\begin{remark}\label{rem1}~
\begin{enumerate}
\item If $(\brond,\|\cdot\|)$ is of stable type $p$, then there exists $q>p$, such that $(\brond,\|\cdot\|)$ is also of stable type $q$ \cite{Maureypisier}.
\item If $(\brond,\|\cdot\|)$ is of stable type $p$, it is also of Rademacher type $p$ \cite{Pisier1}.
\item The above definition of stable type is not the classical one (which involves standard stable r.v.), but an equivalent statement (see for instance \cite{Ledoux1}, proposition 9.12) which will be used in our proofs.
\end{enumerate}
\end{remark}

\begin{remark}~\\
Let $(\brond,\|\cdot\|)$ be a Banach space and $1\leq p\leq 2$. If $\brond$ is of (Rademacher) type $p$, there exists a constant $c(p)$ such that for every finite sequence $(\va X_i)$ with $\|\va X_i\|\in L^p$~:
\begin{equation}\label{relation01}
\EM{\left\|\sum\limits_{i=1}^n\va X_i\right\|}^p\leq c(p)\sum\limits_{i=1}^n\EM\|\va X_i\|^p.
\end{equation}
\end{remark}

The stable type $p$ has been characterized in terms of Marcinkiewicz-Zygmund like SLLN~:

\begin{theorem}[Maurey-Pisier \cite{Maureypisier}]\label{theo2}~\\
Let $1<p<2$. The following two properties are equivalent~:
\begin{enumerate}
\item $(\brond,\|\cdot\|)$ is of stable type $p$.
\item For every bounded sequence $(x_i)$ in $\brond$, the sequence $\left(\frac{1}{n^{1/p}}\sum\limits_{k=1}^n\varepsilon_k x_k\right)$ converges a.s. to 0.
\end{enumerate}
\end{theorem}

\begin{theorem}[Woyczynski \cite{Woyczynski}]\label{theo3}~\\
Let $1<p<2$. The following two properties are equivalent~:
\begin{enumerate}
\item $(\brond,\|\cdot\|)$ is of stable type $p$.
\item For every sequence $(\va X_i)$ of independent, strongly centered $\brond$-valued r.v. for which there exists a nonnegative r.v. $\xi$ with $\EM\xi^p<+\infty$ such that~:
$$
\exists c>0,\ \forall t>0,\ \forall i\in\NM^*,\ \PM(\|\va X_i\|>t)\leq c\PM(\xi>t),
$$
the sequence $\left(\frac{\va S_n}{n^{1/p}}\right)$ converges a.s. to 0, where, as usual, $\va S_n:=\va X_1+\cdots+\va X_n$.
\end{enumerate}
\end{theorem}

In this paper, we will prove a result in the same spirit as Woyczynski's result, but in which the SLLN behaviour of the sequence $\left(\frac{\va S_n}{n^{1/p}}\right)$ is even a quasimartingale behaviour.

\section{A quasimartingale characterization of spaces of stable type $p$.}

We start this section by defining the quasimartingale behaviour of $\left(\frac{\va S_n}{n^{1/p}}\right)$~:

\begin{definition}\label{defi4}~\\
Let $(\va X_k)$ be a sequence of independent, strongly centered $(\brond,\|\cdot\|)$-valued r.v. Denote $\va S_n:=\va X_1+\cdots+\va X_n$ and $\frond_n:=\sigma(\va X_1,\ldots,\va X_n)$. Let $p\in]1,2[$. The sequence $\left(\frac{\va S_n}{n^{1/p}},\frond_n\right)$, or simplier $\left(\frac{\va S_n}{n^{1/p}}\right)$, is a quasimartingale if~:
\begin{equation}\label{1etoile}
\sum\limits_{n=1}^{+\infty}\EM\left\|\left.\EM\left(\frac{\va S_{n+1}}{{(n+1)}^{1/p}}-\frac{\va S_n}{n^{1/p}}\right|\frond_n\right)\right\|<+\infty.
\end{equation}
\end{definition}

\begin{remark}\label{rem5}~\\
Since the r.v. $(\va X_k)$ are independent and centered, condition (\ref{1etoile}) is equivalent to~:
$$
\sum\limits_{n=1}^{+\infty}\frac{\EM\|\va S_n\|}{n^{1+1/p}}<+\infty.
$$
\end{remark}

Now we are able to state our characterization of spaces of stable type $p$ $(1<p<2)$~:

\begin{theorem}\label{theo6}~\\
Let $1<p<2$. The following two properties are equivalent~:
\begin{enumerate}
\item $\brond$ is a stable type $p$ space;
\item For every sequence $(\va X_n)$ of independent, strongly centered r.v., such that
\begin{equation}\label{2etoiles}
\int_0^{+\infty}g^{1/p}(t)dt<+\infty,
\end{equation}
where
$$
\forall t>0,\ g(t):=\sup\limits_{n\geqslant 1}\PM(\|\va X_n\|>t),
$$
$\left(\frac{\va S_n}{n^{1/p}}\right)$ is a quasimartingale.
\end{enumerate}
\end{theorem}

\begin{remark}\label{rem2}~
\begin{enumerate}
\item Condition (\ref{2etoiles}) is not surprising~: indeed, in the i.i.d. case, it can be written~:
\begin{equation}\label{Etoile}
\int_0^{+\infty}\PM^{1/p}(\|\va X\|>t)dt<+\infty,
\end{equation}
which condition is necessary for $\left(\frac{\va S_n}{n^{1/p}}\right)$ being a quasimartingale in every Banach space $\brond$ (see \cite{hechnerheinkel}, proposition 3).
\item Property (\ref{Etoile}) implies that $\EM\|\va X\|^p<+\infty$ (see \cite{hechnerheinkel}, remark 4).
\item There exists a small class of r.v. $\va X$ such that $\EM\|\va X\|^p<+\infty$ and for which (\ref{Etoile}) does not hold (see example 1 in \cite{hechnerheinkel}).
\end{enumerate}
\end{remark}

So -- comparing theorems \ref{theo3} and \ref{theo6} -- one sees that the price to pay for getting a quasimartingale behaviour for $\left(\frac{\va S_n}{n^{1/p}}\right)$ instead of a simple a.s. convergence to 0, is to sharpen a little bit the hypothesis $\EM\xi^p<+\infty$ of theorem \ref{theo3}.

\bigskip

\proof~

In the sequel, $c_k$ will denote positive constants which precise value does not matter.

Let us show the implication $1\implique 2$.

First consider the special case where there exists $M>0$ such that~: $\forall t>M,\ g(t)=0$. Then $\forall k,\ \|\va X_k\|\leq M$ a.s. The space $\brond$ being of Rademacher type $q$ for some $q>p$ by comments 1) and 2) following definition \ref{defi1}, and using relation (\ref{relation01}), one has~:
\begin{align*}
\frac{\EM\|\va S_n\|}{n^{1+1/p}}&\leq\frac{c(q)}{n^{1+1/p}}{\left(\sum\limits_{k=1}^n\EM\|\va X_k\|^q\right)}^{\frac{1}{q}}\leq\frac{c(q)M n^\frac{1}{q}}{n^{1+\frac{1}{p}}}
\end{align*}
and the series having general terms $\frac{\EM\|\va S_n\|}{n^{1+\frac{1}{p}}}$ converges.

\medskip

From now we suppose that $\forall t>0, g(t)>0$.

By a classical symmetrization argument it suffices to consider the case of symmetrically distributed r.v. $(\va X_k)$. For showing that under condition (\ref{2etoiles}) the series having general term $\frac{\EM\|\va S_n\|}{n^{1+\frac{1}{p}}}$ converges, one will split each r.v. $\va X_1,\ldots,\va X_n$ involved in the sum $\va S_n$ into two parts $\va U_{n,k}$ and $\va V_{n,k}$ by truncating $\va X_k$ at a suitable level $v_n$.

For defining $v_n$, one first notices that, $g$ being decreasing, one has~:
$$
\sup\limits_{t>0}t^pg(t)\leq{\left(\int_0^{+\infty}g^{\frac{1}{p}}(x)dx\right)}^p<+\infty.
$$

Furthermore, multiplying the $\va X_k$ by a suitable constant if necessary, one can suppose without loss of generality that~:
\begin{equation}\label{3etoiles}
\sup\limits_{t>0}t^pg(t)\leq 1.
\end{equation}

Now define $v_n:=\inf\left(t>0|g(t)\leq\frac{1}{n}\right)$.

It follows from the definition of $g$ and (\ref{3etoiles}) that~:
$$
v_n\leq n^{1/p}\qquad\text{and}\qquad g(v_n)\leq\frac{1}{n}.
$$

For every $n\in\NM^*$ and $k=1,\ldots,n$, one considers the following centered r.v. (by symmetry)~:
$$
\va U_{n,k}:=\va X_k\va 1_{(\|\va X_k\|\leq v_n)} \qquad\text{and}\qquad \va V_{n,k}:=\va X_k\va 1_{(\|\va X_k\|> v_n)},
$$
and the associated sums~:
$$
\va A_n:=\sum\limits_{k=1}^n\frac{\va U_{n,k}}{n^{1+\frac{1}{p}}}\qquad\text{and}\qquad\va B_n:=\sum\limits_{k=1}^n\frac{\va V_{n,k}}{n^{1+\frac{1}{p}}}.
$$

For showing that the series having general terms $\frac{\EM\|\va S_n\|}{n^{1+\frac{1}{p}}}$ converges, one will show that~:
\begin{equation}\label{4etoiles}
\sum\limits_{n=1}^{+\infty}\EM\|\va A_n\|<+\infty
\end{equation}
and
\begin{equation}\label{5etoiles}
\sum\limits_{n=1}^{+\infty}\EM\|\va B_n\|<+\infty.
\end{equation}

We first prove (\ref{4etoiles}).

By remark $1.$ following definition \ref{defi1}, there exists $q>p$ such that $(\brond,\|\cdot\|)$ is also $q$--stable.

Now suppose that the r.v. $\va X_k$ are defined on a probability space $(\Omega,\trond,\PM)$ and consider $(\varepsilon_k)$ a sequence of independent Rademacher r.v. defined on another probability space $(\Omega^\prime,\trond^\prime,\PM^\prime)$. By symmetry, one has~:
$$
\EM\|A_n\|=\int_\Omega\left(\int_{\Omega^\prime}\frac{1}{n^{1+\frac{1}{p}}}\left\|\sum\limits_{k=1}^n\varepsilon_k(\omega^\prime)\va U_{n,k}(\omega)\right\|d\PM^\prime(\omega^\prime)\right)d\PM(\omega).
$$
By application of the definition of the stable type $q$, one obtains~:
$$
\EM\|A_n\|\leq\frac{C(q)}{n^{1+\frac{1}{p}}}\EM\left(\|(\|\va U_{n,1}\|,\ldots,\|\va U_{n,n}\|)\|_{q,\infty}\right).
$$
For bounding the tails of the weak-$\ell_p$ norm of a sequence of positive, independent r.v., we will use the following classical result due to Marcus Pisier \cite{Marcuspisier}~:
\begin{lemma}\label{lemme6}~\\
For positive valued, independent r.v. $\xi_1,\ldots,\xi_n$, one has~:
\begin{equation}\label{relp6}
\forall q\geq 1,\ \forall u>0,\ \PM\left(\|(\xi_1,\ldots\xi_n)\|_{q,\infty}>u\right)\leq\frac{2e}{u^q}\Delta(\xi_1,\ldots,\xi_n),
\end{equation}
where $\Delta(\xi_1,\ldots,\xi_n)=\sup\limits_{t>0}\left(t^q\sum\limits_{k=1}^n\PM(\xi_k>t)\right)$.
\end{lemma}

For simplicity, denote $\Delta_n$ the quantity $\Delta(\|\va U_{n,1}\|,\ldots,\|\va U_{n,n}\|)$, and notice that by application of lemma \ref{lemme6}~:
\begin{align*}
\EM\left(\|(\|\va U_{n,1}\|,\ldots,\|\va U_{n,n}\|)\|_{q,\infty}\right)
&=\int_0^{+\infty}\PM\left(\|(\|\va U_{n,1}\|,\ldots,\|\va U_{n,n}\|)\|_{q,\infty}>u\right)du\\
&\leq\Delta_n^{\frac{1}{q}}+\int_{\Delta_n^{\frac{1}{q}}}^{+\infty}\frac{2e}{u^q}\Delta_ndu\\
&\leq\cpp\con\Delta_n^{\frac{1}{q}}
\end{align*}

So the proof of (\ref{4etoiles}) reduces to the following lemma~:

\begin{lemma}\label{lemma7}
$$
\sum\limits_{n=1}^{+\infty}\frac{\Delta_n^{\frac{1}{q}}}{n^{1+\frac{1}{p}}}<+\infty
$$
\end{lemma}

Proof of lemma \ref{lemma7}~:\\

One first notices that~:
\begin{align*}
\Delta_n&\leq\sup\limits_{t\leq v_n}t^q\sum\limits_{k=1}^n\PM(\|\va X_k\|>t)
\leq\sup\limits_{t\leq v_n}nt^qg(t)
\leq n{\left(\int_0^{v_n}g^{\frac{1}{q}}(u)du\right)}^{q},
\end{align*}
the last inequality following from the fact that $g$ is decreasing.

For concluding the proof of lemma \ref{lemma7} it remains to check that the series with general term  $a_n:=\frac{1}{n^{1+\frac{1}{p}-\frac{1}{q}}}\int_0^{v_n}g^\frac{1}{q}(u)du$ converges.

First observe that~:
$$
\sum\limits_{n=1}^{+\infty}a_n\leq\sum\limits_{n=1}^{+\infty}\frac{1}{n^{1+\frac{1}{p}-\frac{1}{q}}}\sum\limits_{j=0}^n\int_{v_j}^{v_{j+1}}g^{\frac{1}{q}}(u)du,
$$
where $v_0:=0$, and then exchange the summations in $n$ and $j$~:
$$
\sum\limits_{n=1}^{+\infty}a_n\leq\cpp\con\left(v_1+\sum\limits_{j=1}^{+\infty}\frac{1}{j^{\frac{1}{p}-\frac{1}{q}}}\left(\int_{v_j}^{v_{j+1}}g^{\frac{1}{q}}(u)du\right)\right),
$$
so, by the definition of $v_{j+1}$~:
$$
\sum\limits_{n=1}^{+\infty}a_n\leq\cpp\con\left(v_1+\int_0^{+\infty}g^{\frac{1}{p}}(u)du\right),
$$
which concludes the proof of lemma \ref{lemma7}.

Now we are going to prove (\ref{5etoiles}).

First notice the following chain of inequalities~:
\begin{align*}
\EM\|\va B_n\|&\leq\sum\limits_{k=1}^n\frac{\EM\|\va V_{n,k}\|}{n^{1+\frac{1}{p}}}
=\sum\limits_{k=1}^n\int_0^{+\infty}\PM\left(\frac{\|\va V_{n,k}\|}{n^{1+\frac{1}{p}}}>u\right)du\\
&\leq\sum\limits_{k=1}^n\frac{v_n}{n^{1+\frac{1}{p}}}\PM(\|\va X_k\|>v_n)+\frac{n}{n^{1+\frac{1}{p}}}\int_{v_n}^{+\infty}g(u)du\\
&\leq\frac{nv_n}{n^{1+\frac{1}{p}}}g(v_n)+\frac{1}{n^{\frac{1}{p}}}\int_{v_n}^{+\infty}g(u)du
\leq\frac{v_n}{n^{1+\frac{1}{p}}}+\frac{1}{n^{\frac{1}{p}}}\int_{v_n}^{+\infty}g(u)du.
\end{align*}

Therefore, for proving (\ref{5etoiles}), it suffices to check that condition (\ref{2etoiles}) implies the convergence of the two series with general terms $\frac{v_n}{n^{1+\frac{1}{p}}}$ and $\frac{1}{n^{\frac{1}{p}}}\int_{v_n}^{+\infty}g(u)du$.

\begin{lemma}\label{lemme8}~\\
If (\ref{2etoiles}) is fulfilled, then $\sum\limits_{n=1}^{+\infty}\frac{v_n}{n^{1+\frac{1}{p}}}<+\infty$.
\end{lemma}

Proof of lemma \ref{lemme8}~:

For $j\in\NM^*$, one denotes $t_j:=\frac{v_j+v_{j+1}}{2}$. Then~:
\begin{align}\label{6etoiles}
\int_0^{+\infty}g^\frac{1}{p}(u)du&\geq\sum\limits_{j=1}^{+\infty}\int_{v_j}^{t_j}g^\frac{1}{p}(t)dt
\geq\sum\limits_{j=1}^{+\infty}\frac{1}{{(j+1)}^{\frac{1}{p}}}\frac{v_{j+1}-v_j}{2}.
\end{align}
Now observe that~:
\begin{align}\label{7etoiles}
\sum\limits_{j=1}^{n-1}\frac{1}{{(j+1)}^{\frac{1}{p}}}(v_{j+1}-v_j)&=-\frac{v_1}{2^{\frac{1}{p}}}+\sum\limits_{j=2}^{n-1}v_j\left(\frac{1}{j^{\frac{1}{p}}}-\frac{1}{{(j+1)}^\frac{1}{p}}\right)+\frac{v_n}{n^\frac{1}{p}}
\end{align}
As~:
$$
\frac{v_n}{2}g^\frac{1}{p}(v_n)\leq\int_{\frac{v_n}{2}}^{v_n}g^{\frac{1}{p}}(u)du,
$$
one gets $\lim\limits_{n\to+\infty} v_ng^{\frac{1}{p}}(v_n)=0$ and also $\lim\limits_{n\to+\infty}\frac{v_n}{n^{\frac{1}{p}}}=0$.

As $\frac{1}{j^\frac{1}{p}}-\frac{1}{{(j+1)}^{\frac{1}{p}}}\geq\frac{\cpp\con}{j^{1+\frac{1}{p}}}$, it follows from (\ref{6etoiles}) and (\ref{7etoiles}) that the series having general term $\frac{v_n}{n^{1+\frac{1}{p}}}$ converges.

For completing the proof of the implication $1\implique 2$ of theorem \ref{theo6}, it remains to check~:

\begin{lemma}\label{lemme9}~\\
Under (\ref{2etoiles}), one has $\sum\limits_{n=1}^{+\infty}\frac{1}{n^{\frac{1}{p}}}\int_{v_n}^{+\infty}g(u)du<+\infty$.
\end{lemma}

Proof of lemma \ref{lemme9}.

Let us write~:

$$
\alpha:=\sum\limits_{n=1}^{+\infty}\frac{1}{n^\frac{1}{p}}\int_{v_n}^{+\infty}g(u)du=\sum\limits_{n=1}^{+\infty}\frac{1}{n^\frac{1}{p}}\sum\limits_{j=n}^{+\infty}\int_{v_j}^{v_{j+1}}g(u)du.
$$
By exchanging the summations in $n$ and $j$, one gets~:
\begin{align*}
\alpha&=\sum\limits_{j=1}^{+\infty}\left(\int_{v_j}^{v_{j+1}}g(u)du\right)\sum\limits_{n=1}^j\frac{1}{n^\frac{1}{p}}
\leq\cpp\con\sum\limits_{j=1}^{+\infty}\left(\int_{v_j}^{v_{j+1}}g^\frac{1}{p}(u)du\right)\frac{j^{1-\frac{1}{p}}}{j^{1-\frac{1}{p}}}\\
&\leq\cpp\con\int_0^{+\infty}g^\frac{1}{p}(u)du.
\end{align*}

Let us now show the converse implication $2\implique 1$ of theorem \ref{theo6}.

We first show a general property which is of independent interest~:

\begin{proposition}\label{prop10}~\\
Let $(\va Y_n)$ be a sequence of independent strongly centered r.v. with values in a general Banach space $(\brond,\|\cdot\|)$. Denote by $\va T_n$ the sum $\va Y_1+\ldots+\va Y_n$ and by $\grond_n$ the $\sigma$-field $\sigma(\va Y_1,\ldots,\va Y_n)$. If $\left(\frac{\va T_n}{n^{\frac{1}{p}}},\grond_n\right)$ is a quasimartingale, then $\left(\frac{\va T_n}{n^{\frac{1}{p}}}\right)$ converges a.s. to 0.
\end{proposition}

Proof of proposition \ref{prop10}~:

As noticed earlier, if $\left(\frac{\va T_n}{n^{\frac{1}{p}}}\right)$ is a quasimartingale, then~:
\begin{equation}\label{rel5}
\sum\limits_{n=1}^{+\infty}\frac{\EM\|\va T_n\|}{n^{1+\frac{1}{p}}}<+\infty.
\end{equation}
By Jensen's inequality~:
$$
\forall N\in\NM^*,\ \sum\limits_{n=N}^{+\infty}\frac{\EM\|\va T_n\|}{n^{1+\frac{1}{p}}}\geq\EM\|\va T_N\|\sum\limits_{n=N}^{+\infty}\frac{1}{n^{1+\frac{1}{p}}}\geq\cpp\con\frac{\EM\|\va T_N\|}{N^\frac{1}{p}},
$$
so by (\ref{rel5}),
\begin{equation}\label{rel6}
\lim\limits_{n\to+\infty}\frac{\EM\|\va T_n\|}{n^{\frac{1}{p}}}=0.
\end{equation}
By the conditionnal version of Jensen's inequality~:
\begin{equation}\label{rel7}
\forall n\in\NM^*,\ \EM\left(\|\va T_{n+1}\||\grond_n\right)\geq\|\va T_n\|,
\end{equation}
so~:
\begin{align*}
\sum\limits_{n=1}^N\EM\left|\EM\left(\left.\frac{\|\va T_{n+1}\|}{{(n+1)}^{\frac{1}{p}}}-\frac{\|\va T_n\|}{n^\frac{1}{p}}\right|\grond_n\right)\right|&\leq\sum\limits_{n=1}^N\EM\left|\EM\left(\left.\frac{\|\va T_{n+1}\|-\|\va T_n\|}{{(n+1)}^{\frac{1}{p}}}\right|\grond_n\right)\right|\\+&\cpp\con\sum\limits_{n=1}^N\frac{\EM\|\va T_n\|}{n^{1+\frac{1}{p}}}
\end{align*}
and by (\ref{rel7})~:
\begin{align*}
\sum\limits_{n=1}^N\EM\left|\EM\left(\left.\frac{\|\va T_{n+1}\|}{{(n+1)}^{\frac{1}{p}}}-\frac{\|\va T_n\|}{n^\frac{1}{p}}\right|\grond_n\right)\right|
&\leq\EM\left(\frac{\|\va T_{N+1}\|}{{(N+1)}^{\frac{1}{p}}}\right)+\cpp\con\sum\limits_{n=1}^N\frac{\EM\|\va T_n\|}{n^{1+\frac{1}{p}}}.
\end{align*}
Finaly, by (\ref{rel5}) and (\ref{rel6}), the sequence $\left(\frac{\|\va T_n\|}{n^{1+\frac{1}{p}}}\right)$ is a positive quasimartingale. Therefore, thanks to theorem 9.4 in \cite{Metivier1}, it converges a.s. to a limit, which, by (\ref{rel6}) is necessary 0.

This concludes the proof of proposition \ref{prop10}.

Let us now come back to the proof of the implication $(2)\implique (1)$ of theorem \ref{theo6}.

Let $(\varepsilon_k)$ be a sequence of independent Rademacher r.v. and $(x_k)$ be a bounded sequence of elements in $\brond$. Defining $M:=\sup\|x_k\|$, $\va X_k:=\varepsilon_k x_k$, one gets~:
$$
\forall t>M,\ g(t)=\sup\limits_{k}\PM(\|\va X_k\|>t)=0,
$$
so condition (\ref{2etoiles}) holds. Therefore $\frac{1}{n^{\frac{1}{p}}}\sum\limits_{k=1}^n\varepsilon_k x_k$ is a quasimartingale, which by proposition \ref{prop10} converges a.s. to 0. The $p$-stability of the space $(\brond,\|\cdot\|)$ then follows from theorem \ref{theo2}.

\section{What happens when $p=1$?}

It is natural to wonder if the spaces of stable type 1 (see \cite{hechner3} for the definition of stable type 1) also admit a \guillegauche quasimartingale characterization \guilledroit. In fact it is the case, by theorem 6 in \cite{hechner3}, which can be reformulated as follows~:

\begin{theorem}
Let $\brond$ be a Banach space. The following two properties are equivalent~:
\begin{enumerate}
\item $\brond$ is of stable type 1.
\item For every sequence $({\bf X}_n)$ of independent, strongly centered r.v., such that
\begin{equation}
\int_0^{+\infty}g(t)dt<+\infty,
\end{equation}
where 
$$\forall t>0,\ g(t):=\sup\limits_{n\geqslant 1}\PM(\|{\bf X}_n\|\ln(1+\|{\bf X}_n\|)>t),$$
$\left(\frac{{\bf S}_n}{n}\right)$ is a quasimartingale.
\end{enumerate}
\end{theorem}

\end{document}